\documentclass{article}

\usepackage{arxiv}

\usepackage[utf8]{inputenc} 
\usepackage[T1]{fontenc}    
\usepackage{hyperref}       
\usepackage{url}            
\usepackage{booktabs}       
\usepackage{amsfonts}       
\usepackage{nicefrac}       
\usepackage{lipsum}

\usepackage{morefloats}
\usepackage{color}
\usepackage{multirow}
\usepackage{latexsym}
\usepackage{amsmath,amssymb,amsthm,graphicx}
\usepackage{dsfont}
\usepackage{enumitem}
\usepackage{hyperref}

\newtheoremstyle{thm}
{9pt}
{9pt}
{\itshape}
{}
{\bfseries}
{.}
{ }
{}
\theoremstyle{thm}

\newtheorem{theorem}{Theorem}[section]

\newtheorem{corollary}[theorem]{Corollary}
\newtheorem{proposition}[theorem]{Proposition}

\newtheoremstyle{def}
{9pt}
{9pt}
{}
{}
{\bfseries}
{.}
{ }
{}
\theoremstyle{def}

\newtheorem{remark}[theorem]{Remark}

\newcommand{\HH}{\mathbb{H}} 

\renewcommand{\footnoterule}{%
	\kern -3.5pt
	\hrule width \textwidth height 1pt
	\kern 3.5pt
}

\makeatletter
\def\blfootnote{\xdef\@thefnmark{}\@footnotetext}
\makeatother

\title{A test for Gaussianity in Hilbert spaces via the empirical characteristic functional}

\author{Norbert Henze\\
Institute of Stochastics, \\
Karlsruhe Institute of Technology (KIT), \\
Englerstr. 2, D-76133 Karlsruhe. \\
\texttt{Norbert.Henze@kit.edu}\\
\And  M. Dolores Jim\'{e}nez--Gamero\\
Department of Statistics and Operations Research\\
University of Seville,\\
Seville, Spain\\
\texttt{dolores@us.es}\\
}

\begin{document}

\date{\today}
\maketitle


\begin{abstract}
Let $X_1,X_2, \ldots$ be independent and identically distributed random elements taking values in
a separable Hilbert space $\mathbb{H}$. With applications for functional data in mind, $\mathbb{H}$ may be regarded as a space of
square-integrable functions,  defined on a compact interval.    We propose and study a novel test of the hypothesis $H_0$ that $X_1$ has some
unspecified non-degenerate Gaussian distribution. The test statistic $T_n=T_n(X_1,\ldots,X_n)$ is based on a measure of deviation between
the empirical characteristic functional of $X_1,\ldots,X_n$ and the characteristic functional of a suitable Gaussian random element of $\mathbb{H}$.
We derive the asymptotic distribution of $T_n$ as $n \to \infty$ under $H_0$ and provide a consistent bootstrap approximation thereof.
Moreover, we obtain an almost sure limit of $T_n$ as well as  a normal limit distribution of $T_n$ under alternatives to Gaussianity.
Simulations show that the new test is competitive with respect to the hitherto few competitors available.
\end{abstract}

\section{Introduction}\label{secintroduc}
The normal distribution continues to play a prominent role, since many statistical procedures
for finite-dimensional data assume an underlying normal distribution. It is thus not surprising that
a myriad of tests  for multivariate normality have been proposed. For some more  recent approaches, see e.g.,  \cite{AR07},
\cite{DH08}, \cite{EB12}, \cite{HJ19}, \cite{HJM19}, \cite{HV19}, \cite{KTO07}, \cite{PU05}, \cite{SR05}, \cite{TH14}, \cite{AE09}, and \cite{VPM16}.
A survey of affine invariant tests for multivariate normality is given in \cite{HE02}.

While some of these (and other) tests make use of certain properties that uniquely determine the normal law,
others are based on a comparison of a nonparametric estimator of a function that characterizes
a probability law with a parametric estimator of that function, obtained under the null hypothesis.
A member of the latter class is the test of Epps and Pulley \cite{EP83}. Although originally designed for the univariate case,
the approach of Epps and Pulley was extended to test for multivariate normality by Baringhaus and Henze \cite{BH88} and Henze and Zirkler \cite{HZ90}.
  The resulting procedure, which is usually referred to as the BHEP test, is
        based on a comparison of  the empirical characteristic function (ecf) associated with suitably standardized data, with the characteristic function (cf) of the standard normal law in $\mathbb{R}^d$.
 Because of its nice properties (see Section \ref{BHEP}), the
    BHEP test has been extended in several directions, such as testing for normality of the errors in linear models \cite{JMP05},
     in nonparametric regression models \cite{HM10}, and in GARCH models \cite{JG14}, just to cite a few.

The assumption of normality is important not only in the so-called {\em classical} context,
 in which the data take values in $\mathbb{R}^d$  for some fixed $d \in \mathbb{N}$, but also  in other settings,
 such as functional data analysis. In fact, there are some inferential procedures, designed for functional data,
 that assume {\em Gaussianity} (which is a synonym for normality in that context). Examples are the test for the equality of covariance operators in Ch. 5 of
\cite{HK12}, or the test in  \cite{ZLX10} for the equality of means. On the other hand, some methods are valid under quite general assumptions,
but they greatly simplify when the assumption of normality is added, see, e.g., \cite{BRS18}.
Thus, the problem of testing for Gaussianity is also of interest when dealing with functional data.

The literature of goodness-of-fit tests in the context of functional data, or more general, of data taking values in infinite-dimensional
separable Hilbert spaces, is still rather sparse.
Cuesta-Albertos et al. \cite{CBFM07} consider a test based on random projections,  while Bugni et al. \cite{BHHN09} study an extension
of the Cram\'er--von Mises test. The test  in  \cite{BHHN09}  assumes that the distribution in the null hypothesis depends on a finite-dimensional parameter.
G\'orecki et al. \cite{Gorecki} propose  Jarque-Bera type tests for Gaussianity. For a simple null hypothesis, Ditzhaus and Gaigall \cite{DG18}
employ a test statistic that integrates, along all possible projections, univariate Cram\'er--von Mises  test statistics, obtained by projecting the data.
Since the probability distribution of a random element taking values in a separable Hilbert space is uniquely determined
by its characteristic functional (see Laha and Rohatgi \cite{LR79}), the objective of this paper is to extend the BHEP test to the functional data context.

The paper is organized as follows: Section \ref{BHEP} reviews the BHEP test. Section \ref{test.statistic}  highlights some
key differences between the finite-dimensional case and the functional data context, and it introduces the test statistic,
 whose almost sure limit  is derived in Section \ref{almost.sure.limit}. The asymptotic null distribution of the
 test statistic is obtained in Section \ref{null}. Since this limit distribution depends on
unknown quantities, we prove the consistency of a suitable bootstrap approximation.
In Section \ref{secalt}, we derive the asymptotic distribution of the test statistic under alternatives.
Section \ref{numbers} presents the results of a simulation study, designed to assess the finite-sample
performance of the test, and to compare it with some competitors. It also shows a real data set application. The paper concludes with some remarks.

Throughout the manuscript, we will make use of the following standard notation: The Euclidean norm in $\mathbb{R}^d$, $d \in \mathbb{N}$, will be denoted by
$\| \cdot \|$. The superscript $\top$ means
transposition of column vectors and matrices. We write
N$_d(\mu,\Sigma)$ for the $d$-variate normal
distribution with mean vector $\mu$ and non-degenerate covariance matrix $\Sigma$, and ${\cal N}_d$ stands for the class
of all non-degenerate $d$-dimensional normal distributions. The symbol I$_d$ denotes the unit matrix of order $d$, and i$=\sqrt{-1}$
is the imaginary unit. All random vectors and random elements will be defined on a sufficiently rich probability space
$(\Omega,{\cal A},\mathbb{P})$. The symbols $\mathbb{E}$ and $\mathbb{V}$ denote expectation and variance, respectively, and $ \overset{\mathcal{D}}{=}$ and
 $\overset{\mathcal{D}}{\to}$ mean equality in distribution and convergence in distribution of random vectors and random elements, respectively.
All limits are taken when  $n \rightarrow \infty$, where $n$ denotes the sample size.

%
%
%
\section{The BHEP test in $\mathbb{R}^d$ revisited} \label{BHEP}
In this section, we revisit the BHEP test for finite-dimensional data and take a further view, that will be useful for our purposes.
To this end, let $X_1, \ldots, X_n, \ldots $ be  independent and identically distributed (iid) copies of a $d$-variate random column vector $X$.
We assume that the distribution $\mathbb{P}^{X}$ of $X$ is absolutely continuous with respect to Lebesgue measure.
For testing the hypothesis $H_{0,d}: \ \mathbb{P}^{X}  \in {\cal{N}}_d$,
the rationale of the BHEP test is as follows: Write $\overline{X}_n = n^{-1}\sum_{j=1}^n X_j$ and $S_n = n^{-1}\sum_{j=1}^n (X_j-\overline{X}_n)(X_j-\overline{X}_n)^\top$
for the sample mean and the sample covariance matrix of $X_1,\ldots,X_n$, respectively, and let
$ Y_{n,j} = S_n^{-1/2}(X_j - \overline{X}_n)$, $j =1,\ldots,n$,
be the so-called {\em scaled residuals} of $X_1,\ldots,X_n$, which provide an empirical standardization of $X_1,\ldots,X_n$.
Here, $S_n^{-1/2}$ denotes the unique symmetric square root of $S_n^{-1}$. If $n \ge d+1$, the matrix $S_n$ is invertible with probability one,
see \cite{EP73}. Since, under $H_{0,d}$ and for large $n$, the distribution of the scaled residuals should be close to
the standard $d$-variate normal distribution N$_d(0,{\rm I}_d)$, it is tempting to compare the ecf
\[
\psi_n (t) = \frac{1}{n}\sum_{j=1}^n \exp({\rm i}t^\top Y_{n,j}), \quad t \in \mathbb{R}^d,
\]
of $Y_{n,1},\ldots,Y_{n,n}$ with $\exp(-\|t\|^2/2)$, which is the cf of the law N$_d(0,{\rm I}_d)$.
The BHEP test rejects $H_{0,d}$ for large values of the weighted $L^2$-statistic
\begin{equation} \label{T}
T_{n,d,\beta} = \int_{\mathbb{R}^d}\left | \psi_n(t)-\exp\left(-\frac{1}{2} \|t\|^2\right)\right|^2w_{d,\beta}(t)\, {\rm d}t,
\end{equation}
where $w_{d,\beta}(t) = (\beta^2 2 \pi)^{-d/2} \exp (-\|t\|^2/(2\beta^2))$ is the probability density function
of the $d$-variate normal distribution N$_d(0,\beta^2 {\rm I}_d)$, and $\beta >0$ is a parameter.
In the univariate case, this statistic has been proposed by Epps and Pulley \cite{EP83}, and the extension to
the case $d \ge 2$ has been studied by Baringhaus and Henze \cite{BH88} for the special case $\beta =1$ and, for
general $\beta$, by Henze and Zirkler \cite{HZ90}  and Henze and Wagner \cite{HW97}. The acronym BHEP, after early developers of the idea, was coined by S. Cs\"org\H{o} \cite{CS89},
who proved that the BHEP test is consistent against each non-normal alternative distribution (without any restriction on $\mathbb{P}^X$).
The test statistic $T_{n,d,\beta}$ may be written as
\begin{eqnarray*}
T_{n,d,\beta} & = & \frac{1}{n^2} \sum_{j,k=1}^n \exp \left(- \frac{\beta^2}{2} \|Y_{n,j}- Y_{n,k}\|^2 \right)\\ \nonumber
& &  - \frac{2}{n(1+\beta^2)^{d/2}} \sum_{j=1}^n \exp\left(- \frac{\beta^2}{2(1+\beta^2)} \|Y_{n,j}\|^2 \right) + \frac{1}{(1+2\beta^2)^{d/2}}.
\end{eqnarray*}
This representation shows that $T_{n,d,\beta}$ is a function of the scalar products $Y_{n,j}^\top Y_{n,k} = (X_j-\overline{X}_n)^\top S_n^{-1}(X_k-\overline{X}_n)$,
$1 \le j,k \le n$, and is thus invariant with respect to full rank affine transformations of $X_1,\ldots,X_n$. Moreover, not even the computation of the
square root $S_n^{-1/2}$ is needed. Affine invariance is a ``soft necessary condition" for any genuine test for normality, since the class
${\cal N}_d$ is closed with respect to such transformations, see \cite{HE02}.

%
%
%

\section{The setting and the test statistic} \label{test.statistic}
Assume that $X$ is  a random element of the separable Hilbert space $\mathbb{H} = L^2([0,1],\mathbb{R})$ of (equivalence classes of) square-integrable real-valued functions,
 defined on the compact interval $[0,1]$, with the inner product $\langle f,g\rangle=\int_0^1 f(t)g(t)\, {\rm d}t$,
 norm $\|f\|_{\mathbb{H}} =\langle f,f\rangle^{1/2}$, $f,g \in \mathbb{H}$, and equipped with the Borel $\sigma$-algebra. Throughout the paper we
  assume that $X$ is  square integrable, i.e., we have $\mathbb{E} \|X \|_{\mathbb{H}}^2 <\infty$. As a consequence, we have $\mathbb{E} \|X \|_{\mathbb{H}} <\infty$,
  and  thus there is a unique mean function $\mu = \mathbb{E}(X) \in \mathbb{H}$, which satisfies  $\mathbb{E} \langle X  , x \rangle =\langle \mu  , x \rangle$
   for each $x\in \mathbb{H}$. It follows that $\mu (t)= \mathbb{E} X (t)$ for almost all $t \in [0,1]$.
  Let $c (s,t)= \mathbb{E}[\{X (s)-\mu (s)\}\{X (t)-\mu (t)\}]$, $s,t \in [0,1]$, stand for the covariance function of $X $,
  and write $\mathcal{C} :\mathbb{H}\mapsto \mathbb{H}$ for the covariance operator of $X $,
  defined as $\mathcal{C}  f=\mathbb{E}\left( \langle X-\mu, f\rangle X \right)$, or equivalently,
  as $\mathcal{C}  f(s)=\int c (s,t)f(t)\, {\rm d}t$ for each $f \in \mathbb{H}$.
  The operator ${\cal C}:\mathbb{H} \to \mathbb{H}$ is  linear, compact, symmetric and positive, and it is of trace class. In what follows,
  we denote this class of operators by ${\cal L}^{+}_{{\rm tr}}(\mathbb{H})$.

  Let $X_1, \ldots, X_n$ be iid copies of $X$. The mean function and the covariance function of $X$ can be consistently estimated by means of
  \[
  \overline{X}_n(t)=\frac{1}{n}\sum_{j=1}^nX_j(t), \quad
   c_n(s,t)=\frac{1}{n}\sum_{j=1}^n\{X_j(s)-\overline{X}_n(s)\}\{X_j(t)-\overline{X}_n(t)\},\]
  $s,t \in[0,1]$, respectively. The sample covariance operator $\mathcal{C}_n$, say, is given by $\mathcal{C}_nf(s)=\int_0^1c_n(s,t)f(t)\, {\rm d}t$, $f\in \mathbb{H}$.

 The characteristic functional of $X$, which uniquely determines the distribution $\mathbb{P}^X$ of $X$,   is defined as the function $\varphi: \mathbb{H} \mapsto \mathbb{C}$ ,
 with $\varphi(f)= \mathbb{E} [\exp \left({\rm i} \langle X,f \rangle \right)]$.
 By definition, $\mathbb{P}^X$ is Gaussian if, and only if, there is a $\mu \in \mathbb{H}$ (the {\em expectation} of $X$)
 and ${\cal C} \in {\cal L}^{+}_{{\rm tr}}(\mathbb{H})$ (the {\em covariance operator} of $X$), such that
 \begin{equation}\label{defcfgauss}
 \varphi(f) = \varphi(f; \mu,\mathcal{C})=\exp \left({\rm i} \langle \mu, f\rangle-\frac{1}{2}\langle \mathcal{C}  f, f \rangle \right), \quad f \in \mathbb{H}.
 \end{equation}
In this case, we write $X  \overset{\mathcal{D}}{=} {\rm N}(\mu,{\cal C})$.

 Based on the data $X_1,\ldots,X_n$, we are interested in  testing
  the hypothesis $H_0$ that $\mathbb{P}^X$ is non-degenerate Gaussian, i.e., in a test of
   \[
  H_0:   \varphi(\cdot)=\varphi(\cdot; \mu,\mathcal{C}),\, \mbox{for some $\mu \in \mathbb{H}$ and some  $\mathcal{C} \in {\cal L}^{+}_{{\rm tr}}(\mathbb{H})$},
  \]
  where $\langle {\cal C}f,f \rangle > 0$ for each $f \in \mathbb{H}$ such that $f \neq 0$.
    %
    Two main problems arise when one tries to extend the BHEP test for functional data.

First, a main difference between the finite-dimensional case and the functional data one is that in the latter case we have strict inclusion
\begin{equation} \label{sp}
{\rm sp}_n  \subsetneq \mathbb{H}, 
\end{equation}
where ${\rm sp}_n  =  {\rm sp}(X_1-\overline{X}_n, \ldots, X_n-\overline{X}_n)$ denotes the set of finite linear combinations
of $X_1- \overline{X}_n, \ldots, X_n-\overline{X}_n$,
while in the $d$-dimensional case ${\rm sp}_n = \mathbb{R}^d$  almost surely for each $n\geq d+1$. This  point has an important implication related to invariance.

Section \ref{BHEP} made the case for affine invariance of any genuine test for normality in $\mathbb{R}^d$.
 In the infinite-dimensional case we have the following: If $X \in \mathbb{H}$ is Gaussian
 with mean $\mu$ and covariance operator $\mathcal{C}$,  and $A: \mathbb{H} \to \mathbb{H}$ is a bounded linear operator,
 then $AX$ is also Gaussian (with mean $A\mu$ and covariance operator $A\mathcal{C}A^*$,  $A^*$ being the adjoint of $A$).
 Therefore, arguing as in the previous section, any genuine test for Gaussianity should be invariant under bounded linear operators.
 However, since \eqref{sp} holds for each fixed $n$, it is not reasonable to impose that the test statistic
 be invariant under any bounded linear operator. 

This lack of invariance entails that the null distribution of any test statistic of $H_0$ depends on the population
parameters $\mu$ and $\mathcal{C}$. Therefore, the critical points must be approximated by (for example) some resampling method.
Since our test statistic (to be defined in a moment) is translation invariant, its distribution does not depend on $\mu$.

Second,  recall that the BHEP statistic in $\mathbb{R}^d$ compares the ecf of the scaled residuals
$Y_{n,1}, \ldots, Y_{n,n}$ with the cf of the standard normal law in $\mathbb{R}^d$. There is, however, no
standard normal law in $ \mathbb{H}$. Nevertheless, the BHEP test statistic \eqref{T} can rewritten in the form
\[ 
T_{n,d,\beta} = \int_{\mathbb{R}^d}\left | \phi_n(t)-\phi(t; \overline{X}_n, S_n)\right|^2 F_\beta({\rm d}t).
\] 
Here, $\phi_n(\cdot)$ stands for the ecf of the data $X_1, \ldots, X_n$, $\phi(\cdot ; \mu, \Sigma)$ is the cf of the distribution N$_d(\mu, \Sigma)$,
and $F_\beta$ is a certain distribution function on $\mathbb{R}^d$, see
Lemma 2 in \cite{CSDA09} for details. In view of the above expression, we consider the test statistic
\begin{equation} \label{direct}
T_n= T_n(X_1,\ldots,X_n) = \int_\HH \left |\varphi_n(f)-\varphi(f;  \overline{X}_n,\mathcal{C}_n) \right|^2 Q({\rm d}f)
\end{equation}
for testing $H_0$. Here, $Q$ is some  suitable probability measure on (the $\sigma$-field of Borel subsets of)
 $ \mathbb{H}$, and $\varphi_n$ is the empirical characteristic functional
 \begin{equation}\label{defphinf}
  \varphi_n(f)=\frac{1}{n}\sum_{j=1}^n \exp({\rm i}\langle f,X_j \rangle ), \quad f \in \HH,
 \end{equation}
 of $X_1,\ldots,X_n$.  Straightforward algebra gives
 \begin{eqnarray}\nonumber
 T_n(X_1,\ldots,X_n)  & = & \frac{1}{n} + \frac{2}{n^2} \sum_{1 \le j < k \le n} \int_\HH \cos\left(\langle f, X_j-X_k\rangle \right) Q({\rm d}f)\\ \label{teststatb}
 & & \quad -  \frac{2}{n} \sum_{j=1}^n \int_\HH \cos \left( \langle f,X_j- \overline{X}_n\rangle \right) \exp\left(- \frac{\langle {\mathcal{C}}_nf,f \rangle}{2} \right) Q({\rm d}f)\\ \nonumber
 & & \quad + \int_\HH \exp\left(-  \langle {\mathcal{C}}_nf,f \rangle \right) Q({\rm d}f).
  \end{eqnarray}
Notice that $T_n$ depends solely on the differences $X_j-X_k$ and $X_j - \overline{X}_n$. Consequently, the 
distribution of $T_n$ does not depend on the unknown
expectation $\mu = \mathbb{E}(X)$ of $X$.

In what follows, we will restrict the probability measure $Q$ to be symmetric with respect to the zero element $0$ of $\mathbb{H}$, i.e., $Q$ is invariant with respect to
the mapping $x \mapsto -x$, $x \in \mathbb{H}$. With this assumption, the addition rule $\cos(\alpha - \beta) = \cos \alpha \, \cos \beta + \sin \alpha \, \sin \beta$
and  considerations of symmetry yield
\begin{equation}\label{reprtno}
n T_n = \int_\mathbb{H} V_n^2(f) \, Q({\rm d}f),
\end{equation}
where
\begin{equation} \label{defvn}
V_n(f) = \frac{1}{\sqrt{n}} \sum_{j=1}^n \Big{\{} \cos \langle f,X_j- \overline{X}_n\rangle +  \sin \langle f,X_j- \overline{X}_n\rangle - \exp\left( - \textstyle{\frac12} \langle {\cal C}_nf,f\rangle \right) \Big{\}}, \quad f \in \mathbb{H}.
\end{equation}
The statistic \eqref{direct} is similar to that considered in \cite{BHHN09},
 which is based on a comparison of the empirical {\em distribution} functional with a parametric estimator
 of that functional, obtained under the null hypothesis. As argued in \cite{BHHN09}, $Q$ must be chosen so that the resulting test statistic is tractable computationally.
To this end, notice that the test statistic $T_n$ is  the expected value of the function $V_n(f)$
with respect to the measure $Q$.
Hence, Monte Carlo integration is an option for computation, provided that $Q$ can be easily sampled
(from a computational point of view).
Thus, if $f_1, \ldots, f_M$ is a random sample from $Q$, for some large $M$,
then $T_n$ can be approximated by
\begin{equation} \label{approx}
T_{L,n} =\frac{1}{M}\sum_{m=1}^M  V_n(f_m). 
\end{equation}


\section{An almost sure limit for $T_n$} \label{almost.sure.limit}
This section deals with an almost  sure limit of $T_n$ under general distributional assumptions.

\begin{theorem} \label{limit.dire}
Let $X_1, \ldots, X_n, \ldots $ be iid copies of a random element $X$ of $\mathbb{H}$ satisfying $\mathbb{E} \|X \|^2_\mathbb{H} <\infty$. Writing
$\varphi_X$  for the characteristic functional of $X$, and letting $\mu$ and  $\mathcal{C}$ denote the expectation and the covariance operator of $X$, respectively, we have
 \begin{equation}\label{asvontauq}
 T_n \stackrel{{\rm a.s.}}{\longrightarrow} \tau_Q =\int_\HH \left |\varphi_X(f)-\varphi(f;  \mu,\mathcal{C}) \right|^2 Q({\rm d}f).
\end{equation}
\end{theorem}

\noindent {\sc Proof.}  Recall $\varphi_n(f)$ from (\ref{defphinf}) and $\varphi(f;\mu,{\cal C})$ from (\ref{defcfgauss}).
To stress the dependence on $\omega \in \Omega$, we write
\begin{equation}\label{inttnbound}
T_n(\omega) = \int_\HH \left |\varphi_n(f,\omega)-\varphi(f;  \overline{X}_n(\omega),\mathcal{C}_n(\omega)) \right|^2 Q({\rm d}f),
\end{equation}
where
\begin{eqnarray*}
\varphi_n(f,\omega) & = & \frac{1}{n} \sum_{j=1}^n \exp \left({\rm i}\langle f,X_j(\omega)\rangle \right),\\
\varphi(f,\overline{X}_n(\omega),{\cal C}_n(\omega)) & = & \exp \left({\rm i} \langle \overline{X}_n(\omega),f\rangle - \textstyle{\frac{1}{2}} \langle {\cal C}_n(\omega)f,f\rangle \right),
\end{eqnarray*}
and $\overline{X}_n(\omega) = n^{-1}\sum_{j=1}^n X_j(\omega)$. Moreover, ${\cal C}_n(\omega)$ is the sample covariance operator based on $X_1(\omega),\ldots,X_n(\omega)$.

Let $D\subset \mathbb{H}$ be a countable dense set. By the strong law of large numbers and the fact that the intersection of a countable collection
of sets of probability one has probability one, there is a measurable subset $\Omega_0$ of $\Omega$ such that $\mathbb{P}(\Omega_0) =1$, and for each $\omega \in \Omega_0$ we have, as $n \to \infty$,
$\overline{X}_n(\omega) \to \mu$, ${\cal C}_n(\omega) \to {\cal C}$, $n^{-1} \sum_{j=1}^n \|X_j(\omega)\|_\mathbb{H} \to \mathbb{E} \|X\|_\mathbb{H}$, and
\begin{equation}\label{theorem1b1}
\lim_{n\to \infty} \varphi_n(g,\omega) = \varphi_X(g) \quad \textrm{for each } g \in D.
\end{equation}
Now, fix $\omega \in \Omega_0$ and $f,g \in \mathbb{H}$, and notice that
\[
|\varphi_n(f,\omega)\! - \! \varphi_X(f)|  \le  |\varphi_n(f,\omega)\! -\!  \varphi_n(g,\omega)| +  |\varphi_n(g,\omega) \!  - \!  \varphi_X(g)|
 + |\varphi_X(g) \! - \! \varphi_X(f)| .
\]
If $g \in D$, \eqref{theorem1b1} and the inequality $|e^{{\rm i}u}-e^{{\rm i}v}| \le |u-v|$, valid for real numbers $u$ and $v$, yield
\begin{eqnarray*}
\limsup_{n\to \infty} |\varphi_n(f,\omega) - \varphi_X(f)| & \le & \limsup_{n\to \infty} |\varphi_n(f,\omega) - \varphi_n(g,\omega)| +  |\varphi_X(g) - \varphi_X(f)|\\
& \le & 2 \, \|f-g\|_\mathbb{H} \,  \mathbb{E}\|X\|_\mathbb{H}.
\end{eqnarray*}
Since $2 \, \|f-g\|_\mathbb{H} \,  \mathbb{E}\|X\|_\mathbb{H}$ can be made arbitrarily small because $D$ is dense, the continuity of the exponential function
entails
\[
\lim_{n\to \infty} |\varphi_n(f,\omega)-\varphi(f;  \overline{X}_n(\omega),\mathcal{C}_n(\omega))|^2 =  |\varphi_X(f)-\varphi(f;  \mu,\mathcal{C})|^2
\]
for each $f \in \mathbb{H}$ if $\omega \in \Omega_0$. Since the integrand in \eqref{inttnbound} is bounded from above by 4, the result follows from
dominated convergence. $\Box$

\vspace*{2mm}
Notice that $\tau_Q$ is nonnegative, and that $\tau_Q$ vanishes under $H_0$.
Since the function $t(\cdot)$ that maps $f$ into $t(f)=\left |\varphi_X(f)-\varphi(f;  \mu,\mathcal{C}) \right|^2$ is  continuous,
 and since $H_0$ does not hold if and only if $t(f)>0$ for some
$f \in \mathbb{H}$, a sufficient condition for the validity of $H_0$ if
 $\tau_Q=0$ is
\begin{equation} \label{condition1Q}
Q(B(f,\varepsilon))>0 \quad \text{for each } f \in  \mathbb{H} \text{ and each } \varepsilon>0,
\end{equation}
where $B(f,\varepsilon)=\{ g \in \mathbb{H}: \|f-g\|_\mathbb{H} \leq \varepsilon\}$. Observe that \eqref{condition1Q} holds if $Q$ is Gaussian.
As a consequence, a reasonable test for Gaussianity should reject $H_0$ for large values of
$T_n$. In this respect, it is indispensable to have some information on the
distribution of $T_n$ under the null hypothesis, or at least an approximation to this distribution.
This will be the topic of the next section.

\section{The limit null distribution of $T_n$} \label{null}
In this section we assume that $H_0$ holds, i.e., $X  \overset{\mathcal{D}}{=} {\rm N}(\mu,{\cal C})$ for some $\mu \in \mathbb{H}$ and some
 ${\cal C} \in {\cal L}^{+}_{{\rm tr}}(\mathbb{H})$. Since the distribution of $T_n$ defined in \eqref{direct} does not depend on $\mu$, we will make the tacit
standing assumption $\mu=0$ in what follows. Let $L^2_Q$ denote the Hilbert space
of (equivalence classes of) measurable functions $\Upsilon:\mathbb{H} \mapsto \mathbb{R}$ satisfying $\int_\mathbb{H} \Upsilon(f)^2  Q({\rm d}f)<\infty$.
The scalar product and the resulting norm in $L^2_Q$ will be denoted by $\langle \Upsilon, \Phi \rangle_Q$
 and $\| \Upsilon \|_Q=\sqrt{\langle \Upsilon, \Upsilon \rangle_Q}$, respectively.
Notice that $L^2_Q$ is separable since $\mathbb{H}$ is separable.

The main result of this section is as follows.

\begin{theorem} \label{limittnnull}
Let $X_1, \ldots, X_n, \ldots $ be iid copies of a Gaussian random element $X$ of $\mathbb{H}$ with covariance operator $\mathcal{C}$.
Assume that  $\int_\mathbb{H} \|f\|_\mathbb{H}^4  \, Q({\rm d}f)<\infty$. Then $nT_n \overset{\mathcal{D}}{\to} \|{\cal V}\|^2_Q$, where ${\cal V}$ is a centred Gaussian random element of
$L^2_Q$ having covariance kernel
\begin{equation} \label{covV}
\mathbb{E}[{\cal V}(f){\cal V}(g)]  =  \exp \left(-\textstyle{\frac12}(\sigma_f^2+\sigma_g^2)\right) \Big{\{} \exp(\sigma_{f,g}) - 1 - \sigma_{f,g}- \textstyle{\frac12} \sigma_{f,g}^2 \Big{\}},
\end{equation}
 where  $\sigma_{f,g}=\langle  {\cal C} f,g \rangle$ and  $\sigma^2_f=\sigma_{f,f}$, $f,g \in \mathbb{H}$.
\end{theorem}

\noindent {\sc Proof.}  From \eqref{reprtno}, we have $nT_n = \|V_n\|^2_Q$, where $V_n$ is given in \eqref{defvn}. The idea is to approximate
the random element $V_n$ of $L^2_Q$ by a random element $V_{n,0}$  such that $\|V_n- V_{n,0}\|_Q = o_\mathbb{P}(1)$, and $V_{n,0}$ takes the form
\begin{equation}\label{approxvn}
V_{n,0}(f) = \frac{1}{\sqrt{n}} \sum_{j=1}^n \Psi(f,X_j),
\end{equation}
where $\Psi: \mathbb{H} \times \mathbb{H} \to \mathbb{R}$ is some measurable function satisfying  $\mathbb{E} \Psi(f,X) = 0$ for each $f \in \mathbb{H}$ and
\begin{equation}\label{t2momc}
\mathbb{E} \|\Psi(\cdot, X_1)\|_Q^2 = \int_\mathbb{H} \mathbb{E}\big{[}\Psi^2(f,X_1)\big{]} \, Q({\rm d}f) < \infty.
\end{equation}

In the sequel, the notation $W_n = o_\mathbb{P}(1)$ always refers to  a random element of $L^2_Q$ such that $\|W_n\|_Q$ tends to zero
in probability as $n \to \infty$. Starting with \eqref{defvn}, the addition theorems for the cosine and the sine function yield
\begin{eqnarray}\label{cosexpa}
\cos \langle f,X_j-\overline{X}_n \rangle  & = &  \cos \langle f,X_j \rangle   + \langle f,\overline{X}_n \rangle \, \sin  \langle f,X_j \rangle   + O\left( \langle f,\overline{X}_n \rangle^2\right), \\ \label{sinexpa}
\sin \langle f,X_j-\overline{X}_n \rangle  & = &  \sin \langle f,X_j \rangle  - \langle f,\overline{X}_n \rangle \, \cos \langle f,X_j \rangle   + O\left( \langle f,\overline{X}_n \rangle^2\right).
\end{eqnarray}
Moreover, we have
\begin{equation}\label{diffexpcn}
\exp\left( -\textstyle{\frac12} \langle {\cal C}_n f,f\rangle \right) = \exp\left(-\textstyle{\frac12}\sigma_f^2\right) \left( 1- \textstyle{\frac{1}{2}}\langle ({\cal C}_n-{\cal C})f,f\rangle \right) + O\left(
 \langle ({\cal C}_n-{\cal C})f,f\rangle ^2\right),
\end{equation}
and it follows that
\begin{eqnarray}\nonumber
V_n(f) & = &  \frac{1}{\sqrt{n}}\sum_{j=1}^n \Big{\{} \cos \langle f,X_j\rangle  +  \sin \langle f,X_j\rangle  - \exp\left(- \textstyle{\frac12} \sigma_f^2\right)\Big{\}}\\ \label{termnullt2}
& &  \quad + \, \langle f,\overline{X}_n \rangle \, \frac{1}{\sqrt{n}} \sum_{j=1}^n \Big{\{} \sin \langle f,X_j\rangle  -  \cos \langle f,X_j\rangle  \Big{\}} \\ \label{termnullt3}
& & \quad + \, \frac12 \, \langle \sqrt{n} ({\cal C}_n-{\cal C})f,f\rangle \exp \left( - \textstyle{\frac12} \sigma_f^2 \right) \ + \ o_\mathbb{P}(1).
\end{eqnarray}
Now, the term figuring in \eqref{termnullt2} equals
$
- n^{-1/2} \sum_{j=1}^n \exp\left(- \textstyle{\frac12 }\sigma_f^2 \right) \langle f,X_j \rangle  + o_\mathbb{P}(1)$.
As for the term figuring in \eqref{termnullt3}, we have
\begin{equation}\label{reprcnff}
\langle \sqrt{n} ({\cal C}_n-{\cal C})f,f\rangle = \frac{1}{\sqrt{n}} \sum_{j=1}^n \left( \langle X_j,f\rangle ^2 - \sigma_f^2 \right) + \sqrt{n} \langle \overline{X}_n, f\rangle^2.
\end{equation}
Upon combining we obtain \eqref{approxvn}, where
\[
\Psi(f,x) =  \cos  \langle f,x\rangle  +  \sin \langle f,x\rangle   - \exp\left(- \textstyle{\frac12} \sigma_f^2 \right) \big{\{} 1 + \langle f,x\rangle
- \textstyle{ \frac12} \left( \langle f,x\rangle ^2 - \sigma_f^2 \right)\big{\}}
\]
and $\|V_n - V_{n,0}\|_Q = o_\mathbb{P}(1)$.
From  the central limit theorem in separable Hilbert spaces (see, e.g., Theorem 2.7 in  Bosq \cite{Bosq}), there is a centred
Gaussian random element ${\cal V}$ of $L^2_Q$ with covariance kernel $\mathbb{E}[{\cal V}(f) {\cal V}(g)] = \mathbb{E}[\Psi(f,X)\Psi(g,X)]$,
such that  $V_{n,0} \overset{\mathcal{D}}{\to} {\cal V}$ in $L^2_Q$. From Sluzki's lemma, we thus have
$V_{n} \overset{\mathcal{D}}{\to} {\cal V}$ in $L^2_Q$, and the assertion follows from the continuous mapping theorem. Using the fact that
the joint distribution of $\langle f,X\rangle $ and $\langle g,X\rangle $ is the joint distribution of
$\sigma_f N_1$ and $\sigma_g(\rho N_1+ \sqrt{1-\rho^2}N_2)$, where $N_1,N_2$ are independent standard normal
random variables and $\rho = \sigma_{f,g}/(\sigma_f \sigma_g)$, one easily obtains
\begin{eqnarray*}
\mathbb{E}[\langle f,X\rangle \sin  \langle g,X\rangle ] & = & \sigma_{f,g} \exp\left(- \textstyle{\frac12} \sigma_g^2 \right),\\
\mathbb{E}[\cos \langle f,X\rangle  \langle g,X\rangle^2] & = & (\sigma_g^2 - \sigma_{f,g}^2 ) \exp\left(- \textstyle{\frac12} \sigma_f^2 \right),\\
\mathbb{E}[ \langle f,X\rangle^2 \langle g,X\rangle^2] & = & \sigma_f^2 \sigma_g^2 + 2 \sigma_{f,g}^2,
\end{eqnarray*}
and straightforward algebra shows that $\mathbb{E} [\Psi(f,X)\Psi(g,X)] = \mathbb{E}[{\cal V}(f){\cal V}(g)]$.
Notice that the condition $\int_\mathbb{H} \|f\|_\mathbb{H}^4  \, Q({\rm d}f)<\infty$ implies the validity of \eqref{t2momc}. $\Box$

\vspace*{2mm}

\begin{remark}
Notice that the kernel given in  \eqref{covV} is in accordance with the kernel figuring in display (2.3) of \cite{HW97},
if one replaces $\sigma_f^2$ with $ \|s\|^2$, $\sigma_g^2$ with $\|t\|^2$ and  $\sigma_{f,g}$ with $s^\top t$.
\end{remark}

\vspace*{2mm}

  Since the asymptotic null distribution of  $nT_{n}$  depends on the unknown covariance ope\-ra\-tor $\mathcal{C}$,
     it cannot be used to approximate the actual null distribution of $nT_n$.  To this end, we consider a parametric bootstrap estimator,
defined as follows:
Given $X_1, \ldots, X_n$, let $ {X}^*_{1},  \ldots, {X}^*_{n}$ be iid copies of  $X^*  \overset{\mathcal{D}}{=} {\rm N}(0,{\cal C}_n)$.
Let $T_{n}^*$ be the bootstrap version of $T_n$, which is  obtained by replacing $ {X}_{ 1},  \ldots,  {X}_{n}$
with $ {X}^*_{ 1},  \ldots,  {X}^*_{n}$  in the expression of $T_n$ given in \eqref{teststatb}. Let $\mathbb{P}_*$ denote the conditional distribution, given $X_1,\ldots,X_n$,
and let $\mathbb{P}_0$ denote the null distribution. The bootstrap  estimates
$\mathbb{P}_0(nT_n \leq t)$  by means of $\mathbb{P}_*(nT_n ^*\leq t)$. The next result gives the limit law of the  bootstrap distribution of $nT_n$.

\begin{theorem} \label{limittboot}
Let $X_1, \ldots, X_n, \ldots $ be iid copies of a random element $X$ of $\mathbb{H}$ with covariance operator $\mathcal{C}$. 
Assume that  $\int_\mathbb{H} \|f\|_\mathbb{H}^4  \, Q({\rm d}f)<\infty$. Then $nT_n^* \overset{\mathcal{D}}{\to} \|{\cal V}\|^2_Q$ $ \mathbb{P}^X$-almost surely,
where ${\cal V}$ is the  centred Gaussian random element given in the statement of Theorem \ref{limittnnull}.
\end{theorem}

\noindent {\sc Proof.} Let $V_n^*$ be the bootstrap version of $V_n$ in \eqref{defvn}, which is  obtained by replacing
$ {X}_{ 1},  \ldots,  {X}_{n}$, $\overline{X}_n $ and ${\cal C}_n$ with $ {X}^*_{ 1},  \ldots,  {X}^*_{n}$,
$\overline{X}_n^*$   and ${\cal C}_n^*$, respectively, where $\overline{X}_n^*$  is the sample mean   and ${\cal C}_n^*$ denotes
the sample covariance operator associated with the bootstrap sample  $ {X}^*_{ 1},  \ldots,  {X}^*_{n}$.
 Then $nT_n^*=\|V_n^*\|_Q^2$. Proceeding as in the proof of Theorem \ref{limittnnull}, one obtains
$\|V_n^*- V_{n,0}^*\|_Q = o_{\mathbb{P}_*}(1)$ $\mathbb{P}^X$-almost surely, where 
\[ 
V_{n,0}^*(f) = \frac{1}{\sqrt{n}} \sum_{j=1}^n \widehat{\Psi}_n(f,X_j^*),
\] 
\[
\widehat{\Psi}_n(f,x) =
 \cos  \langle f,x\rangle  +  \sin \langle f,x\rangle   - \exp\left(- \textstyle{\frac12} \widehat{\sigma}_{n,f}^2 \right) \big{\{} 1 + \langle f,x\rangle
- \textstyle{ \frac12} \left( \langle f,x\rangle ^2 - \widehat{\sigma}_{n,f}^2 \right)\big{\}},
\]
and $\widehat{\sigma}_{n,f}^2= \langle {\cal C}_n f, f \rangle $. It thus only remains to show that
$V_{n,0}^* \overset{\mathcal{D}}{\to} {\cal V}$ in $L^2_Q$  $\mathbb{P}$-a.s. With this aim, we apply Theorem 1.1 in Kundu et al. \cite{Kundu}.
To verify the conditions (i)--(iii) of that theorem, let ${\cal C}_{n,V}^*$ and  $c_{n,V}^*$ be the covariance operator and the covariance kernel of $V_{n,0}^*$, respectively.
Likewise, let ${\cal C}_{\cal V}$ and  $c_{\cal V}$ denote the  covariance operator and the covariance kernel of ${\cal V}$, respectively. Notice that
$c_{n,V}^*$ has the same expression as  $c_{\cal V}$ in \eqref{covV}, with $\sigma_f^2$, $\sigma_g^2$ and $\sigma_{f,g}$
replaced by $\widehat{\sigma}_{n,f}^2$, $\widehat{\sigma}_{n,g}^2$ and $\widehat{\sigma}_{n,f,g}=\langle {\cal C}_n f, g \rangle$, respectively.
Notice also that
$c_{n,V}^*(f,g) \stackrel{{\rm a.s.}}{\longrightarrow} c_{\cal V}(f,g)$, for each $f,g \in \mathbb{H}$,
and that $c_V^*(f,g)$ is a bounded function, i.e., for some finite constant $M$ we have $|c_V^*(f,g)| \leq M$ for each $f,g$.
 Let $\{e_k\}_{k \geq 1}$ be an orthonormal basis of $L^2_Q$. By dominated convergence, it follows that
\begin{eqnarray*}
\lim_{n\to \infty} \langle {\cal C}_{n,V}^* e_k,e_\ell\rangle_Q & = & \lim_{n \to \infty}  \int_\mathbb{H} c_{n,V}^*(f,g)e_k(f)e_\ell(g)\, Q({\rm d}f)Q({\rm d}g)\\
       & = & \int_\mathbb{H}  c_{\cal V}(f,g)e_k(f)e_\ell(g) \, Q({\rm d}f) Q({\rm d}g)=\langle{\cal C}_{\cal V} e_k,e_\ell\rangle_Q\quad \mathbb{P}{\rm -a.s.}
\end{eqnarray*}
Setting $a_{k,\ell}=\langle {\cal C}_{\cal V}e_k,e_\ell\rangle_Q$ in the notation of Theorem 1.1 of \cite{Kundu}, this proves
that condition (i) holds. Let $ \mathbb{E}_*$ denote the conditional  expectation, given $X_1,\ldots,X_n$.
 To verify condition (ii) of Theorem 1.1 of \cite{Kundu}, we use Beppo Levi's theorem, Parseval's relation and dominated convergence and obtain
\begin{eqnarray*}
\lim_{n\to \infty} \sum_{k \geq 1} \langle{\cal C}_{n,V}^* e_k,e_k\rangle_Q & = & \lim_{n\to \infty} \sum_{k \geq 1}
\mathbb{E}_* \langle V_{n,0}^*,e_k\rangle^2_Q =\lim_{n \to \infty} \mathbb{E}_*  \| V_{n,0}^*\|_Q^2 \\
 & = & \int_\mathbb{H}  \lim_{n\to \infty}   c_{n,V}^*(f,f) \, Q({\rm d}f) =\int_\mathbb{H}  c_{\cal V}(f,f) \, Q({\rm d}f) = \mathbb{E} \| {\cal V} \|_Q^2
<\infty
\end{eqnarray*}
$\mathbb{P}^X$-almost surely. Finally, we must prove that $L_n(\varepsilon,\Theta) \to 0$ for each $\varepsilon>0$ and each $\Theta \in L^2_Q$, where
\begin{eqnarray*}
L_n(\varepsilon,\Theta) & = & \sum_{j=1}^n \mathbb{E}_* \Big{[} \langle \textstyle{\frac{1}{\sqrt{n}}}\widehat{\Psi}_n( \cdot,X_j^*), \Theta \rangle_Q^2
{\bf 1}\big{\{} \big{|} \langle \textstyle{\frac{1}{\sqrt{n}}} \widehat{\Psi}_n( \cdot,X_j^*), \Theta \rangle_Q \big{|}>\varepsilon\big{\}} \Big{]}\\
& = & \mathbb{E}_* \Big{[} \langle \widehat{\Psi}_n( \cdot,X_1^*), \Theta \rangle_Q^2
{\bf 1} \big{\{} \big{|}\langle  \widehat{\Psi}_n( \cdot,X_1^*), \Theta \rangle_Q\big{|}>\varepsilon \sqrt{n}\big{\}} \Big{]},
\end{eqnarray*}
and  ${\bf 1}\{ \cdot \}$  stands for the indicator function. In the sequel, let $\Theta \neq 0$ without loss of generality.
Using the inequality $t \exp(-t/2) \le 2/e$, $t \ge 0$, it follows that $|\widehat{\Psi}_n(f,x)| \le 4 + \|f\|_\mathbb{H} \|x\|_\mathbb{H} + \textstyle{\frac12} \|f\|^2_\mathbb{H}\|x\|^2_\mathbb{H}$
and thus $\widehat{\Psi}_n^2(f,x) \le \sum_{j=0}^4 a_j \|f\|_\mathbb{H}^j \, \|x\|_\mathbb{H}^j$ for each $f,x \in \mathbb{H}$,  where $a_0,\ldots,a_4$ are positive constants.
By the Cauchy--Schwarz inequality, we have $\langle \widehat{\Psi}_n( \cdot,X_1^*), \Theta \rangle_Q^2 \le \|\widehat{\Psi}_n( \cdot,X_1^*)\|^2_Q  \|\Theta \|^2_Q$. Moreover,
$\big{|}\langle  \widehat{\Psi}_n( \cdot,X_1^*), \Theta \rangle_Q\big{|}>\varepsilon \sqrt{n}$ implies $\|\widehat{\Psi}_n( \cdot,X_1^*)\|^2_Q > \varepsilon^2 n/\|\Theta\|^2_Q$
and thus $\sum_{j=0}^4 a_j \|X_1^*\|^j_\mathbb{H} \int_\mathbb{H}\|f\|^j_\mathbb{H} Q({\rm d}f) > \varepsilon^2 n/\|\Theta\|^2_Q$. As a consequence, we have
\[
{\bf 1} \big{\{} \big{|}\langle  \widehat{\Psi}_n( \cdot,X_1^*), \Theta \rangle_Q\big{|}>\varepsilon \sqrt{n}\big{\}} \le \sum_{k=0}^4 {\bf 1}\Big{\{}
\|X_1^*\|^k_\mathbb{H} > \frac{n\varepsilon^2}{5 a_kb_k\|\Theta\|^2_Q} \Big{\}},
\]
where $b_k = \int_\mathbb{H}\|f\|^k_\mathbb{H} Q({\rm d}f)$, $k \in \{0,\ldots,4\}$, and thus $L_n(\varepsilon,\Theta) \to 0$  will follow if
 \begin{equation}\label{bootts}
 \lim_{n\to \infty} \mathbb{E}_* \big{[} \|X_1^*\|^j_\mathbb{H} {\bf 1}\big{\{} \|X_1^*\|^k_\mathbb{H} >  cn \big{\}}\big{]} = 0, \qquad j,k\in \{0,\ldots,4\},
\end{equation}
where $c$ is a positive constant. Now, \eqref{bootts} holds trivially if $k=0$, and if $k >0$ we have to show
that $\mathbb{E}_* \big{[} \|X_1^*\|^j_\mathbb{H} {\bf 1}\big{\{} \|X_1^*\|_\mathbb{H} >  (cn)^{1/k} \big{\}}\big{]}$ tends to zero $\mathbb{P}^X$-almost surely
for each $j\in \{0,\ldots,4\}$.
The latter convergence follows since $X_1^*  \overset{\mathcal{D}}{=} {\rm N}(0,{\cal C}_n)$  and ${\cal C}_n \to {\cal C}$ $\mathbb{P}^X$-almost surely.
Hence,  condition (iii) of Theorem 1.1 of \cite{Kundu} holds and thus $V_{n,0} \overset{\mathcal{D}}{\to} {\cal V}$ in $L^2_Q$ $\mathbb{P}^X$-almost surely.
 Now, reasoning as in the proof of Theorem \ref{limittnnull}, the result follows. $\Box$

\vspace*{2mm}

Theorem  \ref{limittboot} (which holds regardless of whether $H_0$ is true or not) states that the conditional
distribution of $nT^*_n$ given $X_1,\ldots,X_n$ and the distribution of $nT_n$ when the sample is drawn from a Gaussian population
 with covariance operator ${\cal C}$, are  close to each other for large $n$.
 In particular, under the null hypothesis $H_0$, the conditional distribution of $nT^*_n$ given the data is close to the null distribution of $nT_n$.
More precisely, letting $nT_{n,obs} = nT_n(X_1,\ldots,X_n)$ denote the observed value of the test statistic and, for given $\alpha \in (0,1)$, writing
  $t^*_{n,\alpha}$ for the  upper $\alpha$-percentile of the bootstrap distribution of $nT_n$, the test 
function
$$\Psi_n^*=\left\{
         \begin{array}{ll}
           1, & \text{if}\quad nT_{n,obs} \geq t^*_{n,\alpha}, \\
           0, & \text{otherwise},
         \end{array}
       \right.$$
 or, equivalently, the test that rejects $H_0$ when
$\mathbb{P}_*\{nT^*_n\geq nT_{n,obs}\}\leq \alpha$,
is asymptotically correct in the sense that when $H_0$ is true, we have
$\lim_{n\to \infty} \mathbb{ P}(\Psi_n^*=1) = \alpha$.

An immediate consequence of Theorems \ref{limit.dire} and \ref{limittboot}  is that, if $Q$ satisfies \eqref{condition1Q},
then the test $\Psi_n^*$ is consistent, i.e., it is able to detect any fixed alternative in the sense that
$\lim_{n\to \infty} \mathbb{P}(\Psi_n^*=1) = 1$ whenever $X$ is not Gaussian.

\vspace*{2mm}

 In practice, the  bootstrap distribution of $nT_n$ cannot be calculated exactly. It can  be approximated, however, as follows:
\begin{enumerate} \itemsep=0pt
\item Generate a bootstrap sample $X_1^*, \ldots, X_n^*$, where $X_1^*, \ldots, X_n^*$ are iid from ${\rm N}(0,{\cal C}_n)$.
 \item Calculate the sample mean $\overline{X}_n^*$ and the sample covariance operator ${\cal C}_n^*$ of
 $X_1^*, \ldots, X_n^*$, and compute $nT_{n}^*=nT_{n}(X_1^*, \ldots, X_n^*)$ as given in \eqref{teststatb}, with $\overline{X}_n$ replaced by
 $\overline{X}_n^*$ and ${\cal C}_n$ replaced by  ${\cal C}_n^*$, respectively.
\item Repeat steps 1--2 $B$ times (say),  thus  obtaining $nT_{n}^{*1}, \ldots, nT_{n}^{*B}$. Approximate the  upper $\alpha$-percentile
of the null distribution of  $nT_n$ by the  upper $\alpha$-percentile of the empirical distribution of $nT_{n}^{*1}, \ldots, nT_{n}^{*B}$.
\end{enumerate}


\section{The limit distribution of $T_n$ under alternatives}\label{secalt}
By Theorem \ref{limit.dire}, we have $T_n \stackrel{{\rm a.s.}}{\longrightarrow} \tau_Q$, where $\tau_Q$ is given in \eqref{asvontauq}. In this section,
we will show that, under slightly more restrictive conditions on the underlying distribution, $\sqrt{n}(T_n - \tau_Q)$ has a centred limit normal distribution.
To this end, we first present an alternative representation of $\tau_Q$. Recall the standing assumption that $Q$ is symmetric. We first notice that $\tau_Q$ does not
depend on the expectation $\mu = \mathbb{E}(X)$ of the underlying distribution, since $\varphi_X(f) = \varphi_{X-\mu}(f) \exp({\rm i}\langle f,\mu \rangle )$ and thus
\[
\big{|} \varphi_X(t) - \varphi(f,\mu,{\cal C})\big{|}^2 = \big{|} \varphi_{X-\mu}(f) - \exp\left(-\textstyle{\frac12} \sigma_f^2\right)\big{|}^2,
\]
where $X-\mu$ is centred. Since the covariance operator is invariant with respect to translations, the result follows.

\begin{proposition}\label{proptauq}
We have
\[
\tau_Q = \|z\|^2_Q = \int_\mathbb{H} z^2(f) \, Q({\rm d}f),
\]
where
\begin{equation}\label{defzvonf}
z(f) = \mathbb{E}[\cos \langle f,X\rangle ] + \mathbb{E}[\sin \langle f,X\rangle ] - \exp\left(- \textstyle{\frac12} \sigma_f^2\right), \quad  f \in \mathbb{H}.
\end{equation}
\end{proposition}

 \noindent {\sc Proof.} In view of the discussion above, we assume w.l.o.g. $\mu =0$.
 Using Fubini's theorem and the symmetry of $Q$, \eqref{asvontauq} entails $\tau_Q =  \tau_Q^{(1)} - \tau_Q^{(2)} + \exp(-\sigma_f^2)$, where
 \begin{eqnarray}\nonumber
\tau_Q^{(1)} & = & \int_\mathbb{H} \left( \int_\Omega e^{{\rm i}\langle f,X(\omega)\rangle } \mathbb{P}({\rm d} \omega) \cdot
\int_\Omega e^{-{\rm i}\langle f,X(\omega')\rangle } \mathbb{P}({\rm d} \omega') \right) Q({\rm d}f)\\ \nonumber
& = & \int_\Omega \int_\Omega \left( \int_\mathbb{H} e^{{\rm i} \langle f,X(\omega)-X(\omega')\rangle} Q({\rm d} f) \right) \mathbb{P}({\rm d} \omega) \mathbb{P}({\rm d} \omega')\\ \label{integ3}
& = & \int_\Omega \int_\Omega \left( \int_\mathbb{H} \cos \langle f,X(\omega)-X(\omega')\rangle   Q({\rm d}f) \right) \mathbb{P}({\rm d} \omega) \mathbb{P}({\rm d} \omega')
  \end{eqnarray}
  and
\begin{eqnarray*}
\tau_Q^{(2)} & = & \int_\mathbb{H} \left( \mathbb{E}\big{[} e^{{- \rm i}\langle f,X\rangle}\big{]}e^{-\sigma_f^2/2} +
\mathbb{E}\big{[} e^{{ \rm i}\langle f,X\rangle}\big{]}e^{-\sigma_f^2/2} \right) Q({\rm d}f)\\
& = & 2 \int_\mathbb{H} \mathbb{E}[\cos \langle f,X \rangle ] e^{-\sigma_f^2/2}\, Q({\rm d}f).
\end{eqnarray*}
Using $\cos(\alpha - \beta) = \cos  \alpha \, \cos  \beta +  \sin \alpha   \, \sin  \beta$ and again Fubini's theorem, the expression given in \eqref{integ3}  equals
\[
\int_\mathbb{H} \left(\big{(} \mathbb{E}[\cos \langle f,X\rangle ] \big{)}^2 + \big{(} \mathbb{E}[\sin \langle f,X\rangle ] \big{)}^2  \right) Q({\rm d}f).
\]
Since the symmetry of $Q$ implies $\int_\mathbb{H} \mathbb{E}[\sin \langle f,X \rangle    ] e^{-\sigma_f^2/2} Q({\rm d}f) =0$, the assertion follows. $\Box$

\vspace*{3mm}

We now present the rationale why the limit distribution of $\sqrt{n}(T_n-\tau_Q)$ under fixed alternatives to Gaussianity is a normal distribution.
The main idea is borrowed from \cite{BEH17}, who consider weighted $L^2$-statistics in a more specialized setting. Putting $V_n^\star(\cdot) = V_n(\cdot)/\sqrt{n}$,
where $V_n(\cdot)$ is given in \eqref{defvn}, display  \eqref{reprtno} and Proposition \ref{proptauq} yield
\begin{eqnarray*}
\sqrt{n}\left(T_n-\tau_Q\right) & = & \sqrt{n}\left(\|V_n^\star \|^2_Q - \|z\|^2_Q\right)  \ = \ \sqrt{n} \langle V_n^\star - z,V_n^\star + z\rangle_Q\\
& = & \sqrt{n} \langle V_n^\star - z, 2z + V_n^\star - z\rangle_Q\\
& = & 2  \langle \sqrt{n}( V_n^\star - z), z\rangle_Q + \frac{1}{\sqrt{n}} \|\sqrt{n}(V_n^\star -z)\|^2_Q.
\end{eqnarray*}
Hence, if we can prove the convergence in distribution of $\sqrt{n}( V_n^\star - z)$ in $L^2_Q$ to a centred Gaussian element ${\cal V}^\star$ of $L^2_Q$, then the continuous mapping
theorem and Sluzki's lemma yield $\sqrt{n}(T_n-\tau_Q) \overset{\mathcal{D}}{\to} 2\langle {\cal V}^\star ,z\rangle_Q$, where the distribution of
$2\langle {\cal V}^\star ,z\rangle_Q$ is centred normal with variance $4 \mathbb{E}[\langle {\cal V}^\star ,z\rangle^2_Q]$.

\begin{theorem} \label{thconvnstar}
Assume that  $\int_\mathbb{H} \|f\|_\mathbb{H}^4 \,  Q({\rm d}f)<\infty$, and
let $X_1, \ldots, X_n, \ldots $ be iid copies of a random element $X$ satisfying $\mathbb{E}\|X\|_\mathbb{H}^4 < \infty$.
Let $V_n^\star(\cdot) = V_n(\cdot)/\sqrt{n}$, where $V_n(\cdot)$ is given in \eqref{defvn}. With $z(\cdot)$ defined in
\eqref{defzvonf}, there is a centred Gaussian random element ${\cal V}^\star$ of $L^2_Q$ with covariance kernel
\[
K^\star(f,g) = \mathbb{E}[\xi(f,X)\xi(g,X)], \quad f,g \in \mathbb{H},
\]
where
\begin{eqnarray*}
\xi(f,x) & = & \cos\langle f,x\rangle - \mathbb{E}\cos \langle f,X\rangle +  \sin\langle f,x\rangle - \mathbb{E}\sin \langle f,X\rangle\\
& & \ + \langle f,x\rangle \, \mathbb{E}[\sin \langle f,X\rangle - \cos \langle f,X\rangle ] + \textstyle{\frac12} e^{-\sigma_f^2/2} \left(\langle f,x\rangle^2 -\sigma_f^2\right), \quad x,f \in \mathbb{H},
\end{eqnarray*}
such that $\sqrt{n} \left(V_n^\star - z\right)  \overset{\mathcal{D}}{\to} {\cal V}^\star$.
\end{theorem}

\noindent {\sc Proof.} Notice that
$\sqrt{n}\left(V_n^\star(f) - z(f)\right) = n^{-1/2} \sum_{j=1}^n \left( R_{n,j}(f) + S_{n,j}(f) - T_{n.j}(f)\right)$,
where
\begin{eqnarray*}
R_{n,j}(f) & = & \cos \langle f,X_j-\overline{X}_n \rangle - \mathbb{E}[\cos \langle f,X\rangle ],\\
S_{n,j}(f) & = & \sin \langle f,X_j-\overline{X}_n \rangle - \mathbb{E}[\sin \langle f,X\rangle ],\\
T_{n,j}(f) & = & \exp\left(-\textstyle{\frac12}\langle {\cal C}_nf,f\rangle \right) -  \exp\left(-\textstyle{\frac12}\langle {\cal C}f,f\rangle \right).
\end{eqnarray*}
Using \eqref{cosexpa}, \eqref{sinexpa}, \eqref{diffexpcn} and \eqref{reprcnff}, straightforward calculations yield
\[
\sqrt{n}\left(V_n^\star(f) - z(f)\right) = \frac{1}{\sqrt{n}} \sum_{j=1}^n \xi(f,X_j) \ + \ o_\mathbb{P}(1).
\]
Since $\mathbb{E}[\xi(f,X)] = 0$, $f \in \mathbb{H}$, and since $\mathbb{E}\|\xi(\cdot,X)\|_Q^2 < \infty$ due to the conditions
$\mathbb{E}\|X\|_\mathbb{H}^4 < \infty$ and  $\int_\mathbb{H} \|f\|_\mathbb{H}^4  Q({\rm d}f)<\infty$, the central limit theorem in
 Hilbert spaces and Sluzki's lemma yield the assertion.  $\Box$

\vspace*{3mm}

\begin{corollary} Under the conditions of Theorem \ref{thconvnstar} we have
\[
\sqrt{n}\left( T_n- \tau_Q\right) \overset{\mathcal{D}}{\to} {\rm N}(0,\sigma^2),
\]
where
\[
\sigma^2 = 4 \int_\mathbb{H}\int_\mathbb{H}  K^\star(f,g) z(f)z(g) \, Q({\rm d}f) Q({\rm d}g).
\]
\end{corollary}

\noindent {\sc Proof.} In view of the reasoning preceding Theorem \ref{thconvnstar}, the proof follows from Fubini's theorem, since
\begin{eqnarray*}
\sigma^2 & = & 4\,  \mathbb{E}[\langle {\cal V}^\star ,z\rangle^2_Q]\\
& = & 4 \, \mathbb{E}\bigg{[}\left(\int_\mathbb{H} {\cal V}^\star(f) z(f) \, Q({\rm d}f) \right) \left(\int_\mathbb{H} {\cal V}^\star(g) z(g) \, Q({\rm d}g) \right) \bigg{]} \\
& = & 4 \, \int_\mathbb{H} \int_\mathbb{H} \mathbb{E}[{\cal V}^\star (f){\cal V}^\star(g)] \, z(f)z(g) \, Q({\rm d}f) Q({\rm d}g). \ \Box
\end{eqnarray*}

\begin{remark}
All our results have been stated under the tacit assumption that realizations of $X_1, \ldots, X_n$, i.e., complete trajectories of
 functions, are observable. In practice, these functions are observed at a finite grid of points, and the
  curves $X_1, \ldots, X_n$ are recovered by using nonparametric techniques, such as local linear regression.
 The statistics are then  calculated from $\widehat{X}_1, \ldots, \widehat{X}_n$, which stand for the resulting curve estimators.
 Under suitable assumptions, all previous results remain valid when the test  statistic is calculated from
 $\widehat{X}_1, \ldots, \widehat{X}_n$, see, e.g., Jiang et al. \cite{JHM19}, in particular the comments made after the proof of their Theorem 2.
\end{remark}

\section{Numerical results} \label{numbers}
In this section, we present the results of a simulation study that has been conducted in order to study the finite-sample performance
of the test for Gaussianity based on $T_n$, and to compare the power of this novel test with respect to competing procedures.
All computations have been carried out using programs written in the {\tt R} language (see \cite{R}), with the help of the package {\tt fda.usc}, see  \cite{fda.usc}.

We first studied the performance of the bootstrap approximation to the null
distribution of $T_n$. With this aim, the following experiment was repeated 1000 times:
Independently of each other, we generated $n=50$ realizations  of a standard Wiener process on $[0,1]$ (denoted by W in Table \ref{tabla2}),  and we calculated
$T_{M,n}$  in \eqref{approx}, where $M=1000$ and $Q$ is the Wiener measure on $L^2([0,1])$.
The associated $p$-value was then obtained by generating 200 bootstrap samples.
This setting has been repeated for an Ornstein-Uhlenbeck process  (denoted by OU in Table \ref{tabla2}), and simulations
have also been run for both scenarios  with the sample size $n=100$.

To study the power, we generated samples from
\begin{equation}\label{altmodels123}
Z(t)=A_0+\sqrt{2}\sum_{j=1}^5C_j \cos(2\pi j t)+\sqrt{2}\sum_{j=1}^5S_j \sin(2\pi j t),
\end{equation}
where $A_0$, $C_1, \ldots, C_5$ and $S_1, \ldots, S_5$ are independent random variables, the distributions of which are
displayed in Table \ref{tabla1}.
\begin{table}[h]
\centering
\caption{Description of alternatives.} \label{tabla1}
\begin{tabular}{|l|l|l|}
\hline
Alternative & half-normal & standard normal\\ \hline
Alt1        & $A_0$, $C_1, \ldots, C_5$, $S_1, \ldots, S_5$ &  \\
Alt2        & $A_0$, $C_1, C_2, C_3$, $S_1, S_2, S_3$ & $C_4, C_5$, $S_4, S_5$\\
Alt3        & $A_0$, $C_1$, $S_1$ & $C_2, \ldots, C_5$, $S_2, \ldots, S_5$\\
\hline
\end{tabular}
\end{table}
We also considered the alternatives Alt1', Alt2' and Alt3'. These are the same as the alternatives given in Table \ref{tabla1},
with the exception that the half normal distribution is throughout replaced with an equal mixture of a half normal
distribution and a standard normal distribution. Likewise, the alternatives denoted by
Alt1", Alt2" and Alt3" originate from throughout replacing the half normal distribution with a Laplace distribution (two-sided exponential distribution).

As competitors to the novel test for Gaussianity based on $T_n$, we considered the Jarque-Bera type test
  in  G\'orecki et al. \cite{Gorecki} for iid data (denoted by $JB$ in the table), and the random projection test of
   Cuesta-Albertos et al. \cite{CBFM07} (denoted by $RP$) with 3, 5, 10 and 40 projections.
   Table \ref{tabla2}  reports both the observed empirical level and the empirical power for the nominal levels of significance
    $\alpha=0.05$ and $\alpha = 0.10$.

From Table \ref{tabla2}, we see that the empirical power of the test based on $T_n$ is quite close to the nominal value,
even for the  moderate sample size $n=50$. As for the power, it is not surprising that there is no test
having highest power against all  alternatives considered. For alternatives Alt1, Alt2, Alt3 and Alt3', the test based on $T_n$
outperforms its  competitors, whereas the latter exhibit higher power against the remaining alternatives.
Notice that all alternatives considered belong to the same basic model \eqref{altmodels123}, in which the distribution of some coefficients is switched.
 The power of the $T_n$-test and of the random projection test change softly as the coefficients are switched.
 In most cases, however, the power of the JB test drops as the alternative becomes closer to $H_0$, i.e.,
 as the number of coefficients with normal distribution increases.

\begin{table}
\centering
\caption{Empirical levels and powers for nominal levels $\alpha=0.05, \, 0.10$.} \label{tabla2}
\begin{tabular}{|l|c|c||c|c||c|c|c|c|c|c|c|c|}
\hline
 & \multicolumn{2}{c||}{$T_{n}$} & \multicolumn{2}{c||}{$JB$} &  \multicolumn{8}{c|}{$RP$}\\
 \cline{6-13}
   & \multicolumn{2}{c||}{} &  \multicolumn{2}{c||}{}  & \multicolumn{2}{c|}{3} & \multicolumn{2}{c|}{5} & \multicolumn{2}{c|}{10} & \multicolumn{2}{c|}{40} \\ \hline
$\alpha$     & 0.05  & 0.10  &  0.05 & 0.10  & 0.05  & 0.10    &  0.05  & 0.10   & 0.05   & 0.10    & 0.05 & 0.10 \\ \hline
\multicolumn{13}{|l|}{$n=50$, empirical level}\\ \hline
  W         & .054 & .106 & .041 & .064 &  .067 &  .114 & .059 & .117 & .067 & .127 & .074 & .130  \\
  OU          & .044 & .109 & .037 & .055 &  .063 &  .107 & .069 & .107 & .069 & .124 & .070 & .115  \\ \hline
\multicolumn{13}{|l|}{$n=50$,  empirical  power}\\ \hline
   Alt1        & .605 & .742 & .384 & .452 & .372 & .493 & .388 & .532 &  .432 & .562 & .469 &  .600\\
    Alt2        & .502 & .667 & .154  & .203 & .355 & .459 & .365 & .492 &  .374 & .513 & .416 &  .557 \\
    Alt3        & .518 & .692 & .067 & .095 & .315 & .432 & .307 & .419 &  .317 & .452 & .365 &  .479 \\ \hline
    Alt1'       & .161 & .248 & .426  & .498 & .174 & .266 & .172 & .291 &  .249 & .363 & .269 &  .382 \\
    Alt2'       & .163 & .250 & .169  &  .217 & .145 & .240 & .155 & .240 &  .149 & .245 & .200 &  .303 \\
    Alt3'       & .170 & .272 & .083  & .115 & .152 & .234 & .135 & .223 &  .134 & .224 & .162 &  .245 \\  \hline
Alt1"          & .426  & .538  &  .612 & .662 & .266 & .376 & .295 & .405 & .303 & .421 & .307 & .420  \\
Alt2"          & .471 &.565 & .568  & .634 & .283 & .395 & .326 & .429 & .336 & .453 & .329 & .455 \\
Alt3"          & .468  &.570  & .472 & .516& .319 & .416 & .348 & .475 & .350 & .479 & .349 & .490  \\ \hline
\multicolumn{13}{|l|}{$n=100$,  empirical level}\\ \hline
     W        &  .048     &   .103    & .053 & .072 &  .059 &  .105 & .054 & .104 & .054 & .102 & .049 & .094 \\
  OU          & .050 & .095         & .051 & .077 &  .058 &  .097 & .063 & .114 & .057 & .107 & .047 & .095 \\ \hline
\multicolumn{13}{|l|}{$n=100$,  empirical power}\\ \hline
    Alt1        & .932 & .967 & .766 & .835 & .675 &  .781 & .725 & .829 & .759 & .872 & .811 & .888  \\
    Alt2        & .910 & .956 & .468 & .538 & .638 &  .734 & .658 & .776 & .695 & .828 & .747 & .851  \\
    Alt3        & .888 & .938 &  .102 & .128& .568 &  .678 & .569 & .701 & .603 & .719 & .714 & .812  \\ \hline
    Alt1'       & .210 & .325 & .765 & .824  & .256 &  .384 & .274 & .431 & .434 & .576 & .493 & .648 \\
    Alt2'       & .244 & .369 & .306  & .391 & .249 &  .351 & .233 & .363 & .266 & .382 & .391 & .522 \\
    Alt3'       & .284 & .413 & .150 & .200  & .259 &  .380 & .256 & .362 & .239 & .359 & .292 & .414  \\ \hline
Alt1"          &  .705  & .787    & .885 & .909 &  .461 & .568 & .485 & .619 & .504 & .638 & .507 & .640   \\
Alt2"          & .702  & .797  & .883 & .915 &   .457 & .600 & .498 & .656 & .522 & .655 & .521 & .637 \\
Alt3"          &   .690  & .798  & .829 & .850 & .455 & .575 & .501 & .627 & .533 & .664 & .543 & .666 \\ \hline
\end{tabular}
\end{table}

\begin{figure}
\begin{center}
\begin{tabular}{cc}
\includegraphics[scale=0.55]{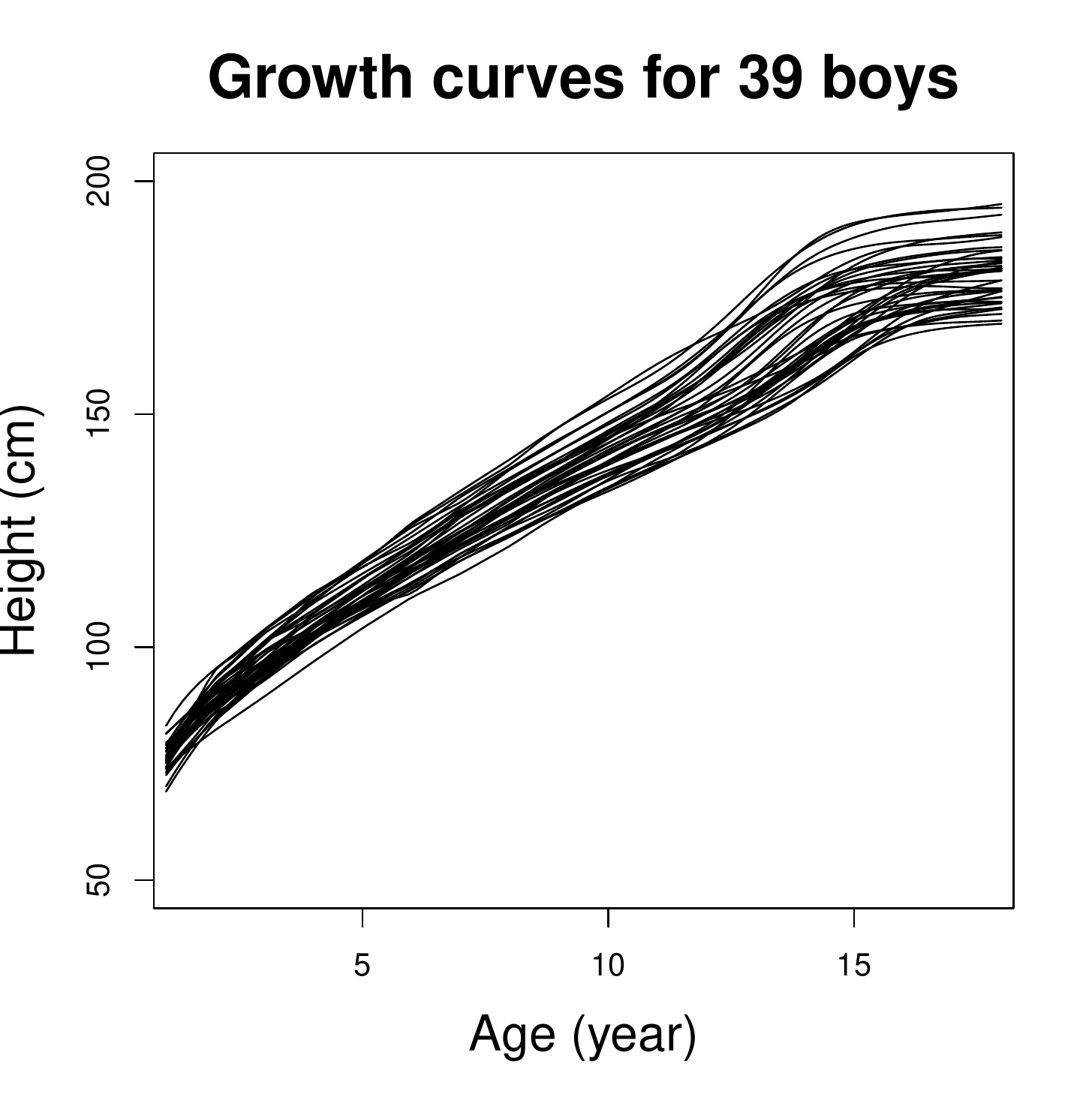} & \includegraphics[scale=0.55]{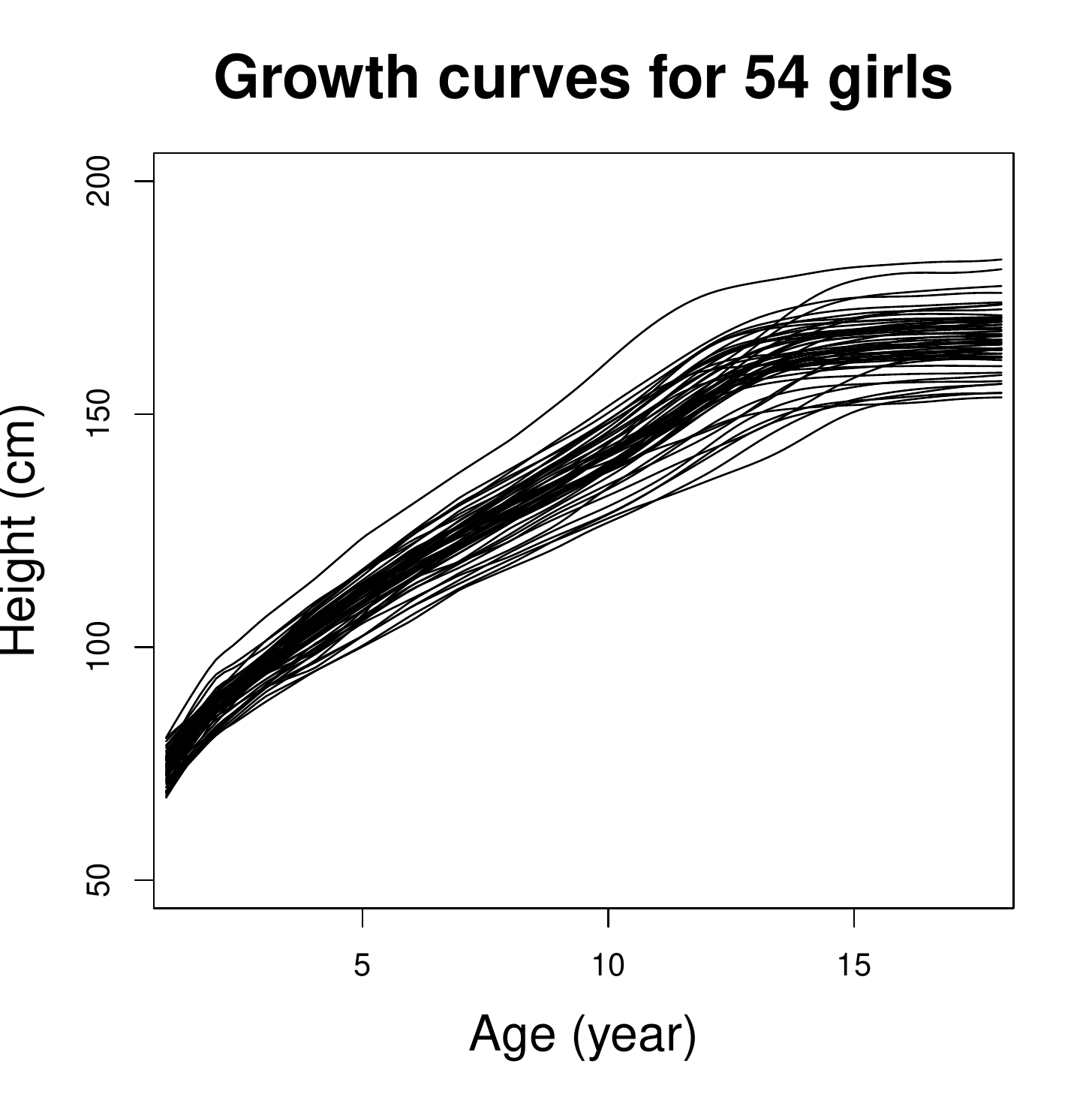}
\end{tabular}
\caption{The Berkeley Growth Data.} \label{curves}
\end{center}
\end{figure}

We close this section with a real data set application. As explained in Section \ref{secintroduc}, some  inferential procedures, designed for
functional data,  assume Gaussianity. An example is the test in Zhang et al. \cite{ZLX10} for the equality of the mean of two functional populations.
Zhang et al. \cite{ZLX10} applied their test to
 the Berkeley Growth Data set. This data set contains the heights of 39 boys and 54 girls recorded at 31 not equally spaced ages from Year 1 to Year 18.
 The data set is available from the {\tt R} package {\tt fda}.
The method in that paper is designed for Gaussian random functions, but the assumption of Gaussianity had not been checked
for either sample.
 We applied the test based on $T_n$ as well as its competitors included in Table \ref{tabla2} to each of the two data sets.
  First of all, proceeding as in \cite{ZLX10}, the growth curves have been
 reconstructed by using local polynomial smoothing. Each of the individual curves has been smoothed separately,
using the same bandwidth $h = 0.3674$. Figure \ref{curves} displays the smoothed growth curves.
Table \ref{p.values} reports the $p$-values obtained. All tests agree in not rejecting the assumption of Gaussianity for both populations, the boys and the girls.

\begin{table}
\centering
\caption{$p$-values for testing Gaussianity for the Berkeley Growth Data.} \label{p.values}
\begin{tabular}{|l|c|c|c|c|c|c|}
\hline
 & \multicolumn{1}{c|}{$T_{n}$} & \multicolumn{1}{c|}{$JB$} &  \multicolumn{4}{c|}{$RP$}\\
 \cline{4-7}
   & \multicolumn{1}{c|}{} &  \multicolumn{1}{c|}{}  & 3 & 5 & 10 & 40 \\ \hline
Boys  &  0.323 & 0.301 & 0.126 & 0.167 & 0.235 & 0.153\\
Girls &  0.144 & 0.433 & 0.107 & 0.138 & 0.153 & 0.248\\ \hline
\end{tabular}
\end{table}
\section{Concluding remarks}\label{secconlus}
We have introduced and studied a novel genuine test for Gaussianity in separable Hilbert spaces that is applicable for functional data.
Some preliminary simulation results show that the procedure compares favorably with the hitherto few existing competitors.
It would be interesting to modify and generalize the approach with respect to testing for Gaussianity  in situations in which the mean and/or the covariance operator
have a certain parametric structure. For example, one could test for a Wiener process on $[0,1]$, where
$\mu =0$ and $c(t,s)=\vartheta \min{t,s}$, for some $\vartheta>0$. It would also be tempting to test for Gaussianity of
 multivariate functional data that  take values in $L^2([0,1]^d)$.

\section*{Acknowledgements}

M.D. Jim\'enez-Gamero has been partially
supported by grant MTM2017-89422-P of the Spanish Ministry of Economy,
Industry and Competitiveness, the State Agency of
Investigation, the European Regional Development Fund.

%
%
%

\begin{small}

\end{small}


\begin{thebibliography}{3}
%
%
%






\bibitem{AR07}
\textsc{Arcones, M.A.}
  (2007).
  Two tests for multivariate normality based on the characteristic function.
\textit{Mathem. Meth. of Statist.},
\textit{16}, 177--201.

\bibitem{BEH17}
\textsc{Baringhaus, L., Ebner, B., and Henze, N.} (2017).
The limit distribution of weighted $L^2$-statistics under alternatives, with applications.
\textit{Ann. Inst. Statist. Math.}, \textit{69}, 969--995.


\bibitem{BH88}
\textsc{Baringhaus, L., and Henze, N.}
  (1988).
  A
consistent test for multivariate normality based on the empirical
characteristic function.
\textit{Metrika},
\textit{35}, 339--348.



\bibitem{BRS18}
\textsc{Boente, G., Rodr\'iguez, D., and Sued, M.}
 (2018).
  Testing equality between several population covariance operators.
  \textit{Ann. Inst. Statist. Math.},
  \textit{70}, 919--950.

\bibitem{Bosq}
\textsc{Bosq, D.}
(2000).
Linear Processes in Function Spaces. Springer, New York.


\bibitem{BHHN09}
\textsc{Bugni, F.A.,  Hall, P.,  Horowitz, J.L., and Neumann, G.R.}
(2009).
Goodness-of-fit tests for functional data.
{\it Econom. J.},
\textit{12}, S1--S18.

\bibitem{CS89}
\textsc{Cs\"org\H{o}, S.}
(1989).
Consistency of some tests for multivariate normality.
\textit{Metrika},
\textit{36}, 107--116.


\bibitem{CBFM07}
\textsc{Cuesta-Albertos, J.A., del-Barrio, E., Fraiman, R., and Matr\'an, C.}
(2007).
The random projection method in goodness of fit for functional data.
\textit{Comput. Statist. Data Anal.},
\textit{51}, 4814--4831.


\bibitem{DG18}
\textsc{Ditzhaus, M., and Gaigall, D.}
(2018),
 A consistent goodness-of-fit test for huge dimensional data and functional data.
 \textit{J. Nonparam. Stat.},
 \textit{30}, 834--859.


\bibitem{DH08}
\textsc{Doornik, J.A., and Hansen, H.}
  (2008).
  An omnibus test for univariate and multivariate normality.
\textit{Oxford Bull. of Economics and Statist.},
\textit{70}, 927--939.



\bibitem{EP73}
\textsc{Eaton, M.L., and Perlman, M.D.}
(1973).
 The non-singularity of generalized sample covariance matrices.
 \textit{Ann.  Statist.},
 \textit{1}, 710--717.


\bibitem{EB12}
\textsc{Ebner, B.}
  (2012).
  Asymptotic theory for the test of multivariate normality by Cox and Small.
\textit{J. Multiv. Anal.},
\textit{111}, 368--379.




\bibitem{EP83}
\textsc{Epps, T.W., and Pulley, L.B.}
(1983).
 A test for normality based on
the empirical characteristic function.
\textit{Biometrika},
\textit{70},
723--726.

%
\bibitem{fda.usc}
\textsc{Febrero-Bande, M., Oviedo de la Fuente, M.} (2012). Statistical Computing in Functional Data Analysis: The R Package fda.usc.
\textit{J. Stat. Softw.} 
\textit{51(4)}, 1--28.

\bibitem{Gorecki}
\textsc{G\'orecki, T., H\"ormann, S., Horv\'ath, L., and Kokoszka P.} (2018)
Testing normality of functional time series.
\textit{J. Time Ser. Anal.}
\textit{39}, 471--487.





\bibitem{HE02}
\textsc{Henze, N.}
(2002).
 Invariant tests for multivariate normality: a critical review.
 \textit{Statist.  Papers},
\textit{43}, 467--506.


\bibitem{HJ19}
\textsc{Henze, N., and Jim\'{e}nez-Gamero, M.D.}
  (2019).
  A class of tests for multinormality with iid and GARCH data based on the empirical moment generating function.
\textit{Test},
\textit{28}, 499--521.


\bibitem{HJM19}
\textsc{Henze, N., Jim\'{e}nez-Gamero, M.D., and Meintanis, S.G.}
  (2019).
  Characterizations of multinormality and corresponding tests of fit, including for GARCH models.
\textit{Econometric Theory},
\textit{35}, 510--546.


\bibitem{HV19}
\textsc{Henze, N., and Visagie, J.}
  (2019).
  Testing for normality in any dimension based on a partial differential equation involving the moment
  generating function.
\textit{Ann. Inst. Statist. Mathem.},
\textit{35}, 510--546. 



\bibitem{HW97}
\textsc{Henze, N., and Wagner, T.}
(1997).
 A new approach to the BHEP tests for multivariate normality.
 \textit{J.  Multiv.  Anal.},
 \textit{62}, 1--23.

\bibitem{HZ90}
\textsc{Henze, N., and Zirkler, B.}
(1990).
A class of affine invariant consistent tests for multivariate normality.
\textit{Commun. Statist. -- Theor. Meth.},
\textit{19}, 3595--3617.



\bibitem{HM10}
\textsc{ Hu\v{s}kov\'a, M., and  Meintanis, S.G.}
(2010).
Test for the error distribution in
nonparametric possiby heteroscedatic regression models.
\textit{Test},
\textit{19}, 92--112.

\bibitem{HK12}
\textsc{Horv\'ath, L., and Kokoszka P.}
(2012).
 Inference for Functional Data with Applications. Springer, New York.

\bibitem{JHM19}
\textsc{Jiang, Q., Hu\v{s}kov\'a, M., Meintanis,  S.G., Zhu, L.}
(2019).
Asymptotics, finite-sample comparisons and applications for two-sample tests with functional data.
\textit{J. Multivar. Anal.},
\textit{170}, 202--220.

\bibitem{JG14}
\textsc{Jim\'enez-Gamero, M.D.}
(2014).
 On the empirical characteristic function process of the residuals in GARCH models and applications.
 \textit{Test},
 \textit{23}, 409--423.

\bibitem{CSDA09}
\textsc{Jim\'enez-Gamero, M.D., Alba-Fern\'andez, M.V.,  Mu\~noz-Garc\'ia, J.,  Chalco-Cano, Y.} (2009)
Goodness-of-fit tests based on empirical characteristic functions.
\textit{Comput. Statist. Data Anal.},
 \textit{53}, 3957--3971.

\bibitem{JMP05}
\textsc{Jim\'enez-Gamero, M.D.,  Mu\~noz-Garc\'ia, J., and Pino-Mej\'ias R.}
(2005).
 Testing
goodness of fit for the distribution of errors in multivariate linear
models.
\textit{J. Multiv. Anal.},
\textit{95}, 301--322.

\bibitem{KTO07}
\textsc{Kankainen, K., Taskinen, S., and Oja, H.}
(2007).
 Tests of multinormality based on location vectors and scatter matrices.
\textit{Statist. Meth. \& Applic.},
\textit{47}, 1923--1934.

\bibitem{Kundu}
\textsc{Kundu, S.,  Majumdar, S. Mukherjee, K.}
(2000).
Central Limit Theorems revisited,
\textit{Statist. Probab. Lett.},
\textit{47}, 265--275.

\bibitem{LR79}
\textsc{Laha, R.G., and Rohatgi, V.K.}
(1979).
 Probability Theory. Wiley, New York.

\bibitem{PU05}
\textsc{Pudelko, J.}
(2005).
 On a new affine invariant and consistent test for multivariate normality.
\textit{Probab. and Mathem. Statist.},
\textit{25}, 43--54.

\bibitem{R}
  R Core Team (2017). R: A language and environment for statistical
  computing. R Foundation for Statistical Computing, Vienna, Austria.
  URL https://www.R-project.org/.



\bibitem{SR05}
\textsc{Sz\'{e}kely, G.J., and Rizzo, M.L.}
(2005).
 A new test for multivariate normality
 \textit{J. Multiv. Anal.},
\textit{93}, 58--80.

\bibitem{TH14}
\textsc{Thulin, M.}
(2014).
 Tests for multivariate normality based on canonical correlations.
\textit{Statist. Meth. \& Applic.},
\textit{23}, 189--208.

\bibitem{AE09}
\textsc{Villase\~nor-Alva, J.A., and Gonz\'alez-Estrada, E.}
  (2009).
  A generalization of Shapiro--Wilk's test for multivariate normality.
\textit{Commun. Statist. -- Theor. Meth.},
\textit{38}, 1870--1883.

\bibitem{VPM16}
\textsc{Voinov, V., Pya, N., Makarov, R., and Voinov, Y.}
(2016).
 New invariant and consistent chi-squared type goodness-of-fit tests for multivariate normality and
 a related comparative study.
\textit{Commun. Statist. -- Theor. Meth.},
\textit{45}, 3249--3263.


\bibitem{ZLX10}
\textsc{Zhang, J.T, Liang, X., and Xiao, S.}
(2010).
 On the two-sample Behrens-Fisher problem for functional data.
 \textit{J. Statist. Theory Pract.},
 \textit{4}, 571--587.


\end{thebibliography}
\end{document}